\documentclass[12pt]{article}

\usepackage{amsfonts}
\usepackage{amssymb,amsmath}

\usepackage[english]{babel}
\usepackage{epsfig}
\usepackage{graphicx}
\catcode`\@=11
\@addtoreset{equation}{section}\catcode`\@=12

\begin{document}

\newtheorem{lm}{Lemma}
\newtheorem{theorem}{Theorem}
\newtheorem{df}{Definition}
\newtheorem{prop}{Proposition}
\newtheorem{rem}{Remark}

\begin{center}
 {\large\bf Birth of discrete Lorenz attractors at the bifurcations of 3D maps
 with homoclinic tangencies to saddle points.}
\vspace{12pt}

 {\bf S.V.Gonchenko}$^1$, {\bf I.I.Ovsyannikov}$^{2}$ and {\bf J.C.Tatjer}$^3$
 \label{Author}
\vspace{6pt}

 {\small
  $^1$ Nizhny Novgorod State University, Russia; \\ E-mail: gonchenko@pochta.ru \\
  $^2$ Universit\"at Bremen, Fachbereich 3, Bibliotekstrasse 1, 
  28359 Bremen, Germany\\
  Nizhny Novgorod State University, Russia;\\
E-mail: ivan.i.ovsyannikov@gmail.com \\
  $^3$ Dept. de Matem\`atica Aplicada i An\`alisi, Universitat de Barcelona,
Spain; \\ E-mail: jcarles@maia.ub.es
 }
\\[12pt]
\end{center}

{\bf Abstract.} It was established in \cite{GMO06} that bifurcations of three-dimensional diffeomorphisms with a homoclinic tangency to
a saddle-focus fixed point with the Jacobian equal to 1 can lead to Lorenz-like strange attractors. In the present paper we prove an
analogous result for three-dimensional diffeomorphisms with a homoclinic tangency to a saddle fixed point with the Jacobian equal to 1,
provided the quadratic homoclinic tangency under consideration is non-simple.


\\~\\

\section*{Introduction}
In this paper we describe a new class of homoclinic tangencies whose bifurcations lead to the birth of strange attractors.
We call the attractor that appears in the Poincare map of a periodically perturbed flow with a Lorenz attractor
a {\em discrete Lorenz attractor}. A theory of such attractors was built in \cite{TS08}. They can emerge in a wide class of maps
which do not need to be directly linked to periodically perturbed flows. Thus, a discrete Lorenz attractor was
first found in \cite{GOST05} for the three-dimensional H\'enon map
\begin{equation}
\bar x = y, \;\; \bar y = z, \;\; \bar z = M_1 + B x + M_2 y - z^2,
\label{HM1}
\end{equation}
where it exists for some open domain of the parameters $(M_1,M_2,B)$ adjoining to
the point $(M_1= - 1/4, B=1, M_2 = 1)$ where  the map has a fixed point with the triplet of multipliers $(-1,-1,+1)$.
It is shown in \cite{SST93} that the normal form for the bifurcations of such fixed point is the Shimizu-Morioka model
subject to exponentially small periodic perturbation; the Shimizu-Morioka model has a region of parameter values which corresponds to
the geometrical Lorenz attractor \cite{ASh86,ASh93}, therefore it is quite typical for a map undergoing the bifurcation
of $(-1,-1,+1)$ to have a discrete Lorenz attractor. Indeed, these attractors were further found numerically in several types of
generalized 3D H\'enon maps \cite{GGOT13}, and in models of non-holonomic mechanics such as Celtic stone \cite{G13}, see also \cite{GG13,GGK13}.
Simple universal bifurcation scenarios that lead to a discrete Lorenz attractor are described in \cite{GGS12,GGKT14}.

Discrete Lorenz attractors belong to the class of the so-called wild hyperbolic attractors \cite{TS98} which
admit homoclinic tangencies and, hence, contain wild hyperbolic sets \cite{N79}, however, bifurcations of these tangencies
do not lead to the birth of periodic sinks. The main reason of this for discrete Lorenz attractors is that they possess
a pseudo-hyperbolic structure. This, very briefly, means that the differential $Df$ of the corresponding map $f$,
in the restriction onto an absorbing neighbourhood ${\cal D}$ of the attractor,
admits an invariant splitting of the form $E^{ss}_x\oplus E^{uc}_x$, depending continuously on the
point $x\in {\cal D}$, such that $Df$ is strongly contracting in restriction to $E^{ss}$ and expands volume
in $E^{uc}$ (see \cite{TS98, TS08} for more detail). This property is robust and prevents from the existence of stable periodic orbits,
therefore the discrete Lorenz attractors preserve  ``strangeness'' at small smooth perturbations. This distinguishes them from
numerous  ``physical'' attractors (quasiattractors in the terminology by Afraimovich and Shilnikov \cite{ASh83a}), such as
H\'enon-like attractors, spiral and screw attractors (R$\ddot{{\rm o}}$ssler attractors, attractors in the
Chua circuits) etc., in which periodic sinks (of arbitrary large periods) can appear under arbitrary small
 perturbations.

The fact that strange (e.g. Lorenz-like) attractors
can appear at the bifurcations of homoclinic tangencies in multidimensional case
was announced yet in \cite{GST93c}. In \cite{GMO06} it was shown that the
discrete Lorenz-like attractors appear at the bifurcations of three-dimensional diffeomorphisms with a homoclinic tangency
to a saddle-focus fixed point with the Jacobian equal to 1. Analogous results were obtained in \cite{GST09,GO10,GO13}
for the bifurcations
of three-dimensional diffeomorphisms with a nontransversal
heteroclinic cycle containing two fixed points, one with the Jacobian less than 1 and the other with the Jacobin greater than 1.

Note that in all these papers it was assumed that at least one of the fixed points is a saddle-focus (i.e. it has a pair
of complex conjugate multipliers inside the unit circle and one real multiplier outside).
This, along with the conditions on the Jacobians of the fixed points, means
that the so-called {\em effective dimension $d_e$ of the problem} (see \cite{T96}) is equal to 3,
which is necessary as the discrete Lorenz-like attractors can exist only in three- and higher-dimensional diffeomorphisms.

In the present paper we study the case of a quadratic homoclinic tangency to a fixed point $O$ of saddle type, i.e.
we assume that all three multipliers of the fixed point are real and different. Then, as it is known from \cite{GST93c,T96,Tat01},
to have the effective dimension of the corresponding problem equal to $3$, we need to assume that
(i) the Jacobian $J$ at the fixed point is equal to $\pm 1$
and (ii) the quadratic homoclinic tangency is {\em non-simple} (see Definition 1 in section~\ref{sec:def}).

The first studies of a codimension-two non-simple homoclinic tangency
were performed in \cite{Tat01} where it was called a {\em generalized} homoclinic tangency.
The notion of a \textit{simple} quadratic homoclinic tangency (a variant of the so-called
quasitransversal homoclinic intersection \cite{NPT}) was introduced in \cite{GST93c}. For three-dimensional maps
with a homoclinic tangency to a saddle fixed point $O$ with multipliers $\nu_i$, $i=1,2,3$ such that $|\nu_1|<|\nu_2|<|\nu_3|$,
the simplicity implies the existence of a non-local two-dimensional invariant manifold, for the map itself and for all $C^1$-close maps.
This manifold contains all orbits entirely lying in a small fixed neighbourhood of the homoclinic orbit.
If the point $O$ has type (2,1), i.e. $|\nu_{1,2}|<1<|\nu_3|$, this manifold is attractive; if the point $O$ has type (1,2),
i.e. $|\nu_{1}|<1<|\nu_{2,3}|$, the manifold is repelling. It follows that neither periodic nor strange attractors
can be born at the bifurcation of a simple tangency if $|\nu_2\nu_3|>1$. However, as it was shown in \cite{Tat01},
if the tangency is non-simple, then periodic attractors can appear provided $|J|=|\nu_1\nu_2\nu_3|<1$, see also  \cite{GGT07}
where the case of a saddle point of type (2,1) was considered in more detail. These results are important for the theory
of dynamical chaos since they show that the non-simple homoclinic tangencies can destroy ``strangeness'' of attractors
as they may lead to the birth of periodic sinks.\\

Let $f_0$ be a three-dimensional orientable $C^r$-diffeomorphism, $r\geq 3$, satisfying the following conditions:

A) $f_0$ has a saddle fixed point $O$ with multipliers $\lambda_1$, $\lambda_2$, $\gamma$
such that $0<|\lambda_2| < |\lambda_1| < 1 < |\gamma|$ (a saddle of type $(2, 1)$);

B) The Jacobian
$J_1 \equiv \lambda_1 \lambda_2 \gamma$ of $f_0$ at the fixed point $O$ is equal to $1$ (a saddle of conservative type);

C) the unstable manifold $W^u(O)$ has a quadratic tangency with $W^s(O)$ at the points of some homoclinic orbit $\Gamma_0$;

D) the tangency is {\em non-simple} (see Definition 1) and nondegenerate.

Diffeomorphisms close to $f_0$ and satisfying  conditions A--D compose, in the space of $C^r$-diffeomorphisms, a locally
connected bifurcation surface of codimension 3. Thus, in order to study bifurcations of $f_0$, we need to consider,
first of all, three-parameter generic unfoldings; the parameters must control the unfolding of degeneracies given by the conditions B, C, and D.
Let $f_{\mu}$ be such a family, where $\mu = (\mu_1,\mu_2,\mu_3)$. We choose the parameters as follows: $\mu_1$ is the splitting parameter
(which controls condition C); $\mu_2$ controls condition D) in such a way that, at $\mu_1=0$,
the tangency becomes simple for $\mu_2 \neq 0$; and $\mu_3$ controls the Jacobian at $O$, i.e it can be taken equal e.g. to
$\mu_3 = 1 - \lambda_1 \lambda_2 \gamma$.
The main result of the present paper is the following\\

{\bf Main Theorem.}  {\em
Let $f_\mu$ be the three-parametric family under consideration {\em (}$f_0$ satisfies  A--D and $f_\mu$ unfolds the degeneracies
given by conditions B, C and D in a generic way{\em)}. Then,
in any neighbourhood of the origin $\mu = 0$ in the parameter space there exist infinitely many domains
$\delta_k$, where $\delta_k\to (0,0,0)$ as $k\to\infty$, such that
the diffeomorphism $f_\mu$  has a discrete Lorenz-like attractor at $\mu\in\delta_k$.} \\

The method of the proof is as follows. For the diffeomorphisms $f_\mu$ we construct first-return maps
$T_k(\mu)$ defined in some neighbourhoods $\sigma^k_0$ near a point of orbit $\Gamma_0$; there are infinitely many such
$\sigma^k_0$ depending on the return time $k$, $k=\bar k, \bar k +1,\dots$. Further, we rescale
coordinates and parameters and show that,
for every sufficiently large $k$, there is an open domain $\Delta_k$ of values of $\mu$ such that (i) $\Delta_k\to (0,0,0)$
as $k\to\infty$ and (ii) the first return map $T_k(\mu)$ for $\mu\in \Delta_k$,
in the rescaled coordinates $(x,y,z)$ and parameters $(M_1,M_2,B)$, takes the form (\ref{HM1})
up to terms that are asymptotically small as $k\to\infty$.
It is important to note that the rescaled coordinates and parameters can take arbitrary finite values (all positive
values for $B$) as $k$ grows. Thus, we can apply the results from \cite{GMO06, GOST05} about the existence of discrete Lorenz-like
attractors in the 3D H\'enon map (\ref{HM1}) and, hence, deduce the existence of such attractor for the map $T_k(\mu)$ for
$\mu\in \delta_k\subset\Delta_k$.

In fact, our analysis provides useful results on global bifurcations in another interesting case of a non-simple
homoclinic tangency. Namely, consider an orientable diffeomorphism $g_0$ which has a saddle fixed point with real multipliers
$\lambda, \gamma_1, \gamma_2$ such that  $0 < |\lambda| < 1 < |\gamma_1| < |\gamma_2|$ and $\lambda \gamma_1 \gamma_2 = 1$
(a conservative saddle of type (1,2)). Suppose also that $g_0$ has a quadratic non-simple
tangency at the points of some homoclinic orbit $\hat\Gamma_0$. Then we can assume that $g_0=f_0^{-1}$ and,
thus, one can obviously use bifurcation results obtained for $f_0$. However, there is an essential difference
in the interpretation of results. Namely, the Main Theorem gives, for $g_0$, only the existence of discrete
Lorenz-like repellers, not attractors.

As $g_\mu$ is inverse to $f_\mu$, the first-return maps $\hat T_k$ for $g_\mu$ are inverse to
the first-return maps $T_k$ for $f_\mu$ and, thus, in the corresponding parameter domains, the rescaled map
$\hat T_k$ is close to the inverse of (\ref{HM1}), i.e. to the map
\begin{equation}
\label{HM2}
\begin{array}{l}
\bar x \; = \; y,\;\; \bar y \; = \; z,\;\;
\bar z = \hat M_1 + \hat M_2 z + \hat B x - y^2,
\end{array}
\end{equation}
where $\displaystyle \hat B = B^{-1}, \hat M_1 = \frac{M_1}{B^2}, \hat M_2 = -\frac{M_2}{B}$.

Map (\ref{HM2}) is well-known in homoclinic dynamics, see e.g. \cite{GST93c, Tat01, GST96, GST08}.
When $\hat{B}=0$ map (\ref{HM2}) becomes two-dimensional: the variable $x$ decouples and we have
(for the coordinates $y$ and $z$) the map of the following form: $\bar y \; = \; z,\;
\bar z = \hat M_1 + \hat M_2 z - y^2$. This map is called Mira map; its dynamics was extensively studied,
see e.g. \cite{Mira96}. In particular, at certain parameter values this non-invertible two-dimensional
map may have a strange attractor with two positive
Lyapunov exponents \cite{PT1,PT2}. Numerical experiments in this case
suggest that the non-robust strange attractors that have the sum of their
Lyapunov exponents positive form a set of large measure in the parameter plane. Numerics also shows that
this property, to have two positive Lyapunov exponents, is inherited by the map (\ref{HM2}) with $B\neq 0$.
The chaotic dynamics of map (\ref{HM2}) was studied e.g. in \cite{GGS12} where it was shown (using the results from \cite{Arneodo79})
that map (\ref{HM2}) possesses strange attractors (quasiattractors)
of spiral type, i.e. those containing a saddle-focus fixed point with two-dimensional unstable manifold.
The question of the existence of genuine strange attractors (e.g. discrete
Lorenz attractors), is open for map (\ref{HM2}), though this problem is very interesting.

The contents of the paper is as follows.
Section~\ref{sec:def} contains the statement of the problem and all necessary
definitions including the definition of two types of non-simple homoclinic
tangencies.
In Section~\ref{res:hom}  we construct the first return maps $T_k$
and formulate the main technical result, Rescaling Lemma~\ref{th3}.
Proof of Lemma~\ref{th3} is given in Section~\ref{sec:th3proof}.

\section{Statement of the problem and main definitions}\label{sec:def}

Let $f_0$ be
a three-dimensional $C^r$-diffeomorphism, $r \geq 3$, satisfying the conditions A)--D). We embed $f_0$ into a  three parameter family $f_\mu$ (general unfolding under conditions B)--D)) with
the parameters $(\mu_1, \mu_2, \mu_3)$ described as above.
We choose a sufficiently small fixed
neighbourhood $U\equiv U(O\cup\Gamma_0)$ of the orbit $\Gamma_0$. Note that $U$  is a union
of a ball $U_0$ containing the point $O$ and a number of
balls surrounding those points of $\Gamma_0$ which lie outside $U_0$.

Denote by  $T_{0}(\mu)$ the restriction of the diffeomorphism $f_\mu$ onto $U_0$.
The map $T_{0}=f_\mu\bigl|_{U_0}$ is called {\it a local map.}
It is known, \cite{GST08, book, GS90, GS92}, that $T_0(\mu)$
can be represented in some $C^r$-smooth local
coordinates $(x_1, x_2, y)$ from $U_0$, smoothly ($C^{r-2}$) depending on $\mu$, in the following {\em main normal form}:
\begin{equation}
\begin{array}{l}
\bar x_1 \; = \; \lambda_1(\mu) x_1 +
\tilde H_1(y, \mu)x_2 + O(\|x\|^2|y|)  \\
\bar x_2 \; = \; \lambda_2(\mu) x_2 + \tilde R_2(x, \mu) +
\tilde H_2(y, \mu)x_2 + O(\|x\|^2|y|) \\
\bar y \; = \; \gamma(\mu) y + O(\|x\||y|^2) ,\\
\end{array}
\label{t0norm}
\end{equation}
where $\tilde H_{1,2}(0,\mu) = 0 \;,\;
\tilde R_{2}(x, \mu) = O(\|x\|^2)$.

We note that, in these coordinates, the local invariant manifolds, stable $W^s_{loc}$,
unstable $W^u_{loc}$ and strong stable $W^{ss}_{loc}$, of the point $O$ are all straightened in $U_0$: their equations are as follows  $W^s_{loc}(O): \; \{ y = 0\}$,
$W^u_{loc}(O): \; \{ x_1 = 0, \; x_2 = 0\}$ and  $W^{ss}_{loc}(O): \; \{ x_1 = 0, \; y = 0\}$.

The intersection points of $\Gamma_{0}$ with $U_0$ belong to
the set $W^{s}\cap W^{u}$ and accumulate at $O$. Thus, infinitely
many points of $\Gamma_0$ lie on $W^{s}_{loc}$ and $W^{u}_{loc}$.
Let $M^{+}(x_1^+,x_2^+,0) \in  W^{s}_{loc}$   and $M^{-}(0,0,y^-) \in W^{u}_{loc}$ be two
of such points and let $M^+=f_0^{n_0}(M^-)$ for some positive
integer $n_0$. Let $\Pi^+\subset U_0$  and $\Pi^-\subset U_0$ be
small enough neighbourhoods of the points $M^{+}$  and $M^{-}$, respectively.
The map $T_{1}(\mu)\equiv f_\mu^{n_0}: \Pi^- \rightarrow  \Pi^+$ is called
{\it a global map.}

In order to formulate explicitly conditions C and D on the homoclinic tangency
we recall first some important facts from the theory of invariant manifolds \cite{book, HPS}.

\begin{figure}
\centerline{\epsfig{file=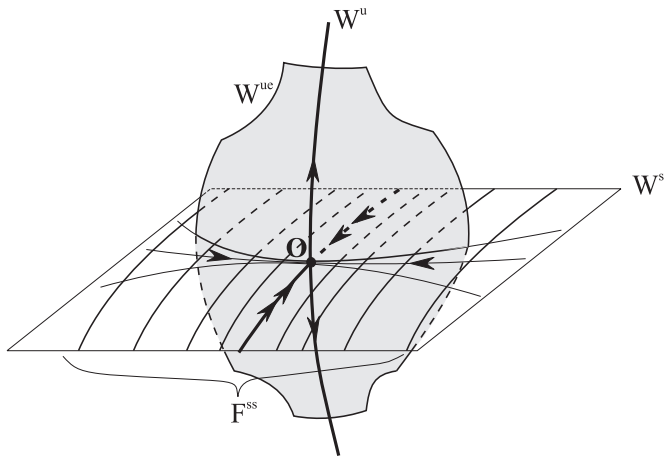,
height=6cm}} \caption{{\footnotesize A part of the strong stable
foliation $F^{ss}$ containing the strong stable manifold $W^{ss}$;
and a piece of one of the extended unstable manifolds $W^{ue}$ containing
$W^u$ and being transversal to $W^{ss}$ at $O$.}} \label{fig04}
\end{figure}

When the condition A is fulfilled, in $U_0$ there exist the so-called
{\em extended unstable manifolds} $W^{ue}(O)$,
see Fig.~\ref{fig04} where an example of such a manifold is presented.
There are infinitely many (continuum) such manifolds, they are
only $C^{1+\varepsilon}$-smooth in general. In the case under consideration, all $W^{ue}(O)$ are
two-dimensional, contain $W^u$ and touch at the point $O$ the eigendirection of $Df_0$ corresponding to the multiplier $\lambda_1$
--- the line $\{y=0,x_2 = 0\}$ in coordinates (\ref{t0norm}). Thus, each
$W^{ue}_{loc}(O)$ has the equation of form $x_2 = \varphi(x_1, y)$, where $\varphi(0, y) \equiv 0$ and
$\varphi'_{x_1}(0, 0) = 0$.
Note that in coordinates (\ref{t0norm})
all the manifolds $W^{ue}(O)$
have the same tangent plane $\{x_2=0\}$ at each point of $W^u(O)$.

Another fact we use
is the existence of a
{\em strong stable invariant foliation}, see Fig.~\ref{fig04}.
Recall that $W^s(O)$ contains the one-dimensional strong stable submanifold $W^{ss}(O)$, which is
invariant, $C^r$-smooth and touches at $O$ the eigenvector corresponding to the strong
stable (nonleading) multiplier $\lambda_2$. Moreover, the manifold $W^s(O)$ is foliated near $O$ by
the foliation $F^{ss}$ which is $C^r$-smooth, unique and contains
$W^{ss}$ as a leaf. As we said before, in  coordinates (\ref{t0norm}), $W^{ss}_{loc}(O)$ has the equation
$\{ x_1=0,\;y = 0\}$ and also the foliation $F^{ss}$ on $U_0$ consists of the leaves
$\{x_1 = \mbox{const},\;y = 0\}$.

Denote the tangent plane to $W^{ue}(O)$ at the point $M^-$ as $P^{ue}(M^-)$.
By our condition, the curve $T_1(W^u_{loc}\cap\Pi^-)$ has at $\mu=0$ a quadratic tangency with $W^s_{loc}$ at the point $M^+$.

\begin{df}
{\rm The homoclinic tangency under consideration is called  \textsf{simple} if $T_{1}(P^{ue}(M^-))$
intersects transversely the leaf $F^{ss}(M^+)$  of the foliation $F^{ss}$ containing the point $M^+$.
If this condition is not fulfilled we call such quadratic tangency \textsf{non-simple}.}
\label{def:simple}
\end{df}

According to \cite{Tat01} we define two general cases of non-simple homoclinic tangencies:

Case ${\rm I}.$ {\it The surface $T_{1}(P^{ue}(M^-))$ is
transversal to the plane $W^{s}_{loc}(O)$  but is tangent to the line
$F^{ss}(M^+)$ at $M^+$.}

Case ${\rm II}.$ {\it The surfaces $T_{1}(P^{ue}(M^-))$ and
$W^{s}_{loc}(O)$ have a tangency at $M_1^+$ and the curves $T_{1}(W^u_{loc}(O) \cap \Pi^-)$
and $F^{ss}(M^+)$ have a general intersection.} \\

\begin{figure}
\centerline{\epsfig{file=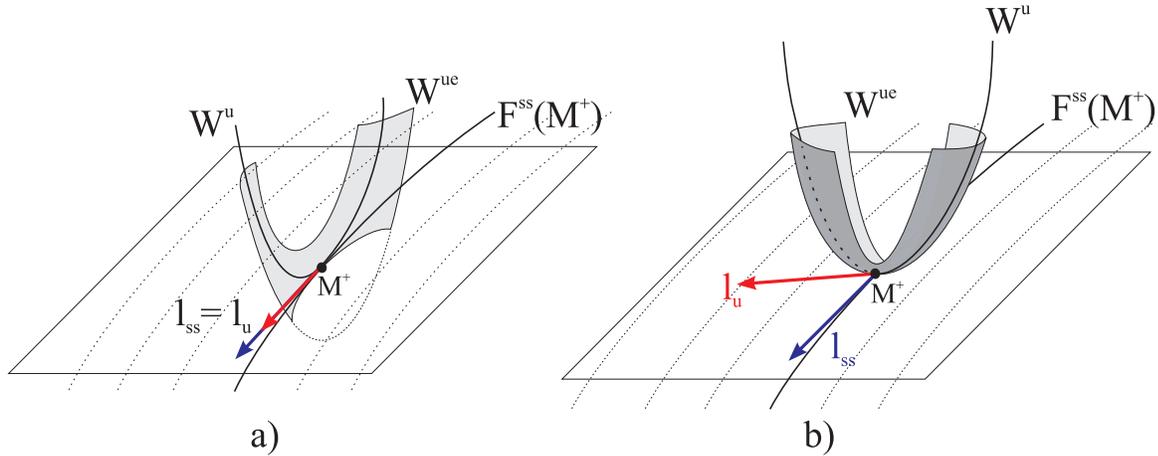,
height=6cm}} \caption{{\footnotesize Two types of nondegenerate non-simple quadratic (homoclinic) tangency:
(a) $W^{ue}$ is transversal to $W^s_{loc}$ and touches the leaf $F^{ss}(M^+)$; (b) $W^{ue}$ is
tangent to $W^s_{loc}$ and the curves $W^u$ and $F^{ss}(M^+)$ has a general intersection at
$M^+$.}} \label{fig05}
\end{figure}

Thus, in Case ${\rm I}$ the tangent vectors $l_u$ to
$T_{1}(W^u_{loc})$ and
$l_{ss}$ to $F_1^{ss}(M_1^+)$ are collinear,
while in Case ${\rm II}$ these vectors have different directions,
see Fig.~\ref{fig05}~(a) and~(b).

\section{Calculation of the first return maps $T_k(\mu)$.
}\label{res:hom}

We consider on $U_0$ the local coordinates $(x_1,x_2,y)$ in which the map $T_0(\mu)$ has the form (\ref{t0norm}). Let $\{p_i(x_{i1},x_{i2},y_{i})\}$, $i=1,2,\ldots,k$, be such points in $U_0$ that $p_{i+1}=T_0(p_i)$. Then, by 
{\rm \cite{GST08, book, GS92}}, we can represent the map $T_0^k: U_0 \to U_0$ in the so-called {\em Shilnikov cross-form} as follows.

\begin{equation}
\begin{array}{l}
x_{k1} - \lambda_1^{k}(\mu) x_{01} =
\hat\lambda^{k}\xi_{k1}(x_0, y_k, \mu)  \\
x_{k2} - \lambda_2^{k}(\mu) x_{02} =
\hat\lambda^{k}\xi_{k2}(x_0, y_k, \mu) \\
y_0 - \gamma(\mu)^{-k} y_k =
\hat\gamma^{-k} \eta_k(x_0, y_k, \mu),
\end{array}
\label{T0kk}
\end{equation}
where $\hat\lambda$ and $\hat\gamma$ are some constants such that
$0 < \hat\lambda < |\lambda_1(\mu)|,\; \hat\gamma >
|\gamma(\mu)|$  and functions $\xi_{k}$ and $\eta_{k}$
are uniformly bounded along with all derivatives up
to order $(r - 2)$. Moreover, $\|x_k\|_{C^{r-1}} = O(\lambda_1^k),\|y_0\|_{C^{r-1}} = O(\gamma^{-k}), \|x_k,y_0\|_{C^{r}} \to 0 $ as $k\to\infty$, see \cite{GGT07} for more details.

To construct the global map $T_1$ we use the facts that
when $\mu = 0$ we have  $T_1(M^-) = M^+$, where, recall, $M^+ = (x_1^+,x_2^+,0)$ and $M^-=(0,0,y^-)$ is a pair of points of $\Gamma_0$, and
$T_1(W_{loc}^u)$ and $W_{loc}^s$ have a quadratic tangency at the point
$M^+$. Accordingly, the global map $T_1$ at all small $\mu$
can be written as
\begin{equation}
\label{t1-1}
\begin{array}{l}
\bar x_1 - x_1^+ \; = \; a_{11}x_1 + a_{12}x_2 + b_1(y - y^-) + O(\|x\|^2+ (y - y^-)^2)
\\
\bar x_2 - x_2^+ \; = \; a_{21}x_1 + a_{22}x_2 +
b_2(y - y^-) + O(\|x\|^2 +(y - y^-)^2) \\
\bar y \; = \; y^+(\mu) + c_1x_1 + c_2x_2 + d(y - y^-)^2 + O(\|x\|^2 + \|x\||y-y^-| + |y - y^-|^3)
\\
\end{array}
\end{equation}
where $y^+(0) = 0$ and
coefficients $a_{11}, \ldots, d$ as well as $x^+$ and $y^-$ depend smoothly on $\mu$. Note that since the homoclinic tangency at $\mu=0$ is quadratic, we have $d \neq 0$. Moreover, map $T_1(0)$ is a
diffeomorphism, therefore
\begin{equation}
\label{t1-3}
J_1 = {\rm det}\; \left(
\begin{array}{rcl}
a_{11}\; & a_{12}\; &\;b_1 \\
a_{21}\; & a_{22}\; &\;b_2 \\
c_1\; & c_2 &\; 0 \\
\end{array}
\right)
\neq 0
\end{equation}
and, hence, $b_1^2 + b_2^2 \neq 0 $, $c_1^2 + c_2^2 \neq 0$ at $\mu=0$.


Now we consider condition D separately for Cases I and II.

Case ${\rm I}$. The tangent plane $P^{ue}(M^-)$ to $W^{ue}_{loc}$ at
point $M^-$ has equation ${x_2 = 0}$. The equation of
$T_1(P_{ue}(M^-))$ at $\mu = 0$ is obtained by putting $x_2 = 0$ into (\ref{t1-1}).
Then the transversality of $T_1(P^{ue}(M^-))$ and $W^s_{loc}$
($\bar y = 0$) yields $c_1(0) \neq 0$. The tangent vector to the
 line $T_1(P_{ue}(M^-)) \cap W^s_{loc}$ at  point $M^+$ 
is $(b_1(0), b_2(0), 0)$. The equation of the leaf $F^{ss}(M^+)$ is $\{x_1 = x_1^+, y = 0\}$.
Therefore, the tangency of $\;T_1(P^{ue}(M^-))$ and $F^{ss}(M^+)$
implies $b_1(0) = 0$. In this case $b_2\neq 0$ and
$a_{11}^2+a_{12}^2\neq 0$ because of (\ref{t1-3}). Thus, in Case I
map $T_1(\mu)$ has form (\ref{t1-1}) where
$$
b_1(0) = 0 \;,\; c_1(0)\neq 0 \;\;{\rm and}\;\; b_2(0)\neq 0.
$$

Case $\rm II$. The equation of $T_1(P^{ue}(M^-))$ at
$\mu = 0$ is the same as in Case~$\rm I$. Then the tangency of
$T_1(P^{ue}(M^-))$ and $W^{s}_{loc}$ at $\mu = 0$ implies
that $c_1(0) = 0$. Also, the tangent vectors to the lines
$T_1(L_u)$ and to $F^{ss}(M^+)$ at point $M^+$ are non-parallel if
$b_1(0) \neq 0$. Thus, in Case~$\rm II$ map $T_1(\mu)$ has form
(\ref{t1-1}) where
$$
c_1(0) = 0 \;,\; b_1(0) \neq 0 \;\;{\rm and}\;\; c_2(0) \neq 0.
$$

The main goal of the paper is to study bifurcations of
single-round periodic orbits of diffeomorphisms close to $f_0$.
Every point of such an orbit can be considered as a fixed point
for the corresponding first return map $T_k$, where $k$ can run all sufficiently large integers. These maps are constructed as the composition $T_k = T_1 T_0^k: \Pi^+ \to \Pi^- \to\Pi^+ $, where, recall, $\Pi^+ \subset U_0$ and $\Pi^- \subset U_0$ are  small neighbourhoods of the homoclinic points $M^+$ and
$M^-$, respectively.
If we take a point $M \in \Pi^+$, then its
iterations under the local map $T_0$ can reach $\Pi^-$. Such
points form on $\Pi^+$ a set consisting of infinitely many
three-dimensional strips $\sigma_k^0 = T_0^{-k}(\Pi^-) \cap \Pi^+$,
$k \in \{k_0, k_0 + 1, \ldots\}$.
Accordingly, $\sigma_k^1 = T_0^k(\sigma_k^0)$ is a
(three-dimensional) strip on $\Pi^-$. Strips $\sigma_k^0$ and
$\sigma_k^1$  accumulate on $W^s_{loc} \cap \Pi^+$ and
$W^u_{loc} \cap \Pi^-$, respectively, as $k \to \infty$. The global
map $T_1$ maps the strip $\sigma_k^1$ into a three-dimensional horseshoe
$T_1(\sigma_k^1) \subset \Pi^+$. By the definition, map $T_k = T_1T_0^k:
\sigma_k^0 \to \Pi^+ \; , \; k \in \{k_0, k_0 + 1, \ldots \},$ is {\em the first
return map}. Since $T_1$ is the $n_0^{th}$ power of $f_\mu$, any fixed
point of $T_k$ corresponds to a single-round periodic orbit of $f_\mu$ of
period $k + n_0$.

We  consider three-parameter families $f_{\mu_1, \mu_2, \mu_3}$ of diffeomorphisms close to $f$.
Naturally, a parameter $\mu_1$ of the splitting of manifolds $W^s(O)$ and $W^u(O)$ with respect to
the point $M^+$ is considered as one of the governing parameters. It is seen from (\ref{t1-1}) that
\begin{equation}
\;\;\mu_1 \equiv y^+(\mu)\;.
\label{mu1}
\end{equation}

As the second parameter, we consider a parameter $\mu_2$ resolving the degeneracy connected either
with condition ${\rm D}_{\rm I}$ (in Case I) or with condition ${\rm D}_{\rm II}$ (in Case II). It
is convenient to take directly
\begin{equation}
\;\mu_2 = b_1(\mu)\;\;\;{\rm in\; Case\; I}
\label{mu2I}
\end{equation}
and
\begin{equation}
\;\mu_2 = c_1(\mu)\;\;\;
{\rm in\; Case\; II.}
\label{mu2II}
\end{equation}

Finally, the third parameter should control the Jacobian $J = \lambda_1\lambda_2\gamma$ of $f_\mu$
at the saddle $O_\mu$. Therefore, we  put for all cases
\begin{equation}
\;\mu_3 = 1 - \lambda_1\lambda_2\gamma
\label{mu3}
\end{equation}

Thus, the family $f_{\mu_1, \mu_2, \mu_3}$ can be considered as a general unfolding of the
corresponding non-simple homoclinic tangency under conditions A--D. \\

 We construct the first return maps $T_k$ using formulae (\ref{T0kk}) and
(\ref{t1-1}). By such a way we  obtain a formula for $T_k$ in the initial (small) variables
$(x_1, x_2, y) \in U_0$ and parameters $\mu_1$ and $\mu_2$. Next, we rescale the initial variables and
parameters
$$
(x_1, x_2, y) \mapsto (X_1, X_2, Y)\;,\; (\mu_1, \mu_2, \mu_3) \mapsto (M_1, M_2, M_3)\;,
$$
with asymptotically small (as $k \to \infty$) factors, in such a way that map $T_k$ is rewritten,  in
the rescaled variables and parameters, as some regular three-dimensional quadratic map (the limit form) plus asymptotically small (as $k \to \infty$) terms.
Moreover, new coordinates $(X_1, X_2, Y)$ and
parameters $(M_1, M_2, M_3)$ can take arbitrary finite values at large $k$.
This result is formulated as the following lemma.

\begin{lm} {\em(}\textsc{Rescaling Lemma}{\em)}.
Let $f_{\mu}$ be the three parameter family under consideration. Then, in the space of the parameters 
there exist infinitely many regions $\Delta_k$
ac\-cu\-mu\-la\-ting to $\mu=0$ as $k \to \infty$, such that the map $T_k$ in appropriate rescaled
coordinates and parameters is asymptotically $C^{r - 1}$-close to one of the following limit maps.

{\rm 1)} In Case I, the limit map is
\begin{equation}
\label{HI}
\begin{array}{l}
\bar X_1 \; = \; - M_3 X_2,\;\; \bar X_2 \; = \; Y,\;\;
\bar Y = M_1 - X_1 + M_2 X_2  - Y^2,
\end{array}
\end{equation}
where
\begin{equation}
\label{M123I}
\begin{array}{l}
M_1 = - d\gamma^{2k}(\mu_1 + \lambda_1^kc_1x_1^+ +\nu_k^1),\; M_2\;=\;  c_1 (\mu_2 +
\rho_k^1)(\lambda_1\gamma)^k, \; M_3 = J_1(\lambda_1\lambda_2\gamma)^k ,
\end{array}
\end{equation}
and $\nu_k^1 = O(|\gamma|^{-k} + \hat\lambda^k)$, $\rho_k^1 = O(|\hat\lambda/\lambda_1|^k)$.

{\rm 2)}  In Case II, the limit map is
\begin{equation}
\label{HII}
\begin{array}{l}
\bar X_1 \; = \; Y,\;\; \bar X_2 \; = \; X_1,\;\;
\bar Y = M_1 + M_2 X_1 + M_3 X_2- Y^2,
\end{array}
\end{equation}
where
\begin{equation}
\label{M123II}
\begin{array}{l}
M_1 = - d\gamma^{2k}[\mu_1 - \gamma^{-k}y^{-} + \lambda_1^k(\mu_2x_1^+ +\nu_k^2)],\; M_2\;=\;  b_1
(\mu_2 + \rho_k^2)(\lambda_1\gamma)^k, \\ M_3 = J_1(\lambda_1\lambda_2\gamma)^k
\end{array}
\end{equation}
and $\nu_k^2,\rho_k^2  = O(|\hat\lambda/\lambda_1|^k + |\lambda_1|^k)$.
\label{th3}
\end{lm}

If $M_3$ is separated from zero,  maps (\ref{HI}) and (\ref{HII}) are equivalent (map (\ref{HI})
takes form (\ref{HII}) after scaling $X_1 \to -M_3 X_1$). Then dynamics of the first return maps for
$\mu \in \Delta_k$ is the same as for the three-dimensional H\'enon map (\ref{HM1}).

Thus, we need only to proof the Rescaling Lemma.

\section{Proof of Lemma~\ref{th3}.}\label{sec:th3proof}

Note that under the
assumptions A and B we have that $|\lambda_1|>|\lambda_2|$, $|\lambda_1\gamma|>1$ and $|\lambda_2\gamma|>1$. In principle, the maps $T_k$ are rescaled differently in the Cases~$\rm I$ and~$\rm II$. However, there is a preparation part
of the proof that is conducted in the same way for both the cases.

\subsection{Preparation form of map $T_k$ for rescaling.}\label{prform}

Using (\ref{t1-1}) and (\ref{T0kk}) one can write the map $T_k =
T_1T_0^k$ for sufficiently large $k$  and small $\mu$ in
the form
\begin{equation}
\label{tk-1}
\begin{array}{l}
\bar x_1 - x_1^+(\mu) \; = \; a_{11}(\lambda_1^{k} x_{1}
+ \hat\lambda^{k}\xi_{k1}(x,y,\mu)) +  a_{12}(\lambda_2^{k} x_{2} + 
\hat\lambda^{k}\xi_{k2}(x,y,\mu)) + \\
 + b_1(y - y^-(\mu)) + 
O((\hat\lambda^{k}+ \lambda_1^{2k})\|x\|^2 + |y - y^-|^2),\\ \\
\bar x_2 - x_2^+(\mu) \; = \;
a_{12}(\lambda_1^{k} x_{1}
+ \hat\lambda^{k}\xi_{k1}(x,y,\mu)) + a_{22}(\lambda_2^{k} x_{2} +
\hat\lambda^{k}\xi_{k2}(x,y,\mu)) + \\
 + b_2(y - y^-(\mu)) + O((\hat\lambda^{k}+ \lambda_1^{2k}) \|x\|^2 + |y - y^-|^2),\\ \\
\gamma^{-k} \bar y  -
\hat\gamma^{-k}\eta_k(\bar x,\bar y,\mu) \; = \; \mu_1 +
c_{1}(\lambda_1^{k} x_{1}
+ \hat\lambda^{k}\xi_{k1}(x,y,\mu)) +  
 c_{2}(\lambda_2^{k} x_{2} + \\
 + \hat\lambda^{k}\xi_{k2}(x,y,\mu))  + d(y - y^-)^2 +
 O((\hat\lambda^{k}+ \lambda_1^{2k})\|x\|^2 + \lambda_1^{k}\|x\||y - y^-| 
  + |y - y^-|^3 )\;.
\\
\end{array}
\end{equation}
We shift coordinates $ x_{1new}\;=\; x_1 - x_1^+(\mu) +
\phi_k^1(\mu),\; x_{2new}\;=\; x_2 - x_2^+(\mu) +
\phi_k^2(\mu),\;  y_{new}\;=\; y - y^-(\mu) +
\psi_k(\mu)$, where $\phi_k,\psi_k = O(|\lambda_1|^{k})$,
in such a way that the right sides of (\ref{tk-1}) do not contain
constant terms for the first two equations and linear in $y_{new}$
terms for the third equation. Then (\ref{tk-1}) takes the form
\begin{equation}
\label{tk-1+}
\begin{array}{l}
\bar x_1 \; = \; a_{11}\lambda_1^{k} x_{1} + b_1 y 
 + (\hat\lambda^{k}+ \lambda_1^{2k})O(\|x\|) + \lambda_1^{k}O(|y|) + O(y^2), \\

\bar x_2  \; = \;
a_{21}\lambda_1^{k}x_{1} + b_2 y +
 (\hat\lambda^{k}+ \lambda_1^{2k})O(\|x\|) + \lambda_1^{k}O(|y|) + O(y^2), \\
\\
\bar y -
(\hat\gamma/\gamma)^{-k}
\eta_k(\bar x + x^+ + \phi_k,\bar y + y^- +
\psi_k,\mu) \; = \;
M_k  + d(1 + s_k^1) \gamma^{k}y^2 + \\
\lambda_1^{k}\gamma^{k}(c_{1} x_{1} +
(|\lambda_1|^{k} + |\hat\lambda/\lambda_1|^{k})O(\|x\|))
+ \lambda_1^k\gamma^k O(\|xy\|) + \gamma^{k} O(y^3), \\
\end{array}
\end{equation}
where $s_k^1 = O(\lambda_1^{k})$ and
\begin{equation}
\label{Mk1}
\begin{array}{l}
M_k\;=\;\gamma^{k}[\mu_1 - \gamma^{-k}(y^- + \dots) + \lambda_1^{k}(c_1x_1^+ +
\dots)]
\end{array}
\end{equation}
and dots stand for coefficients tending to zero as $k\to\infty$.

Consider the third equation of (\ref{tk-1+}).
First of all, we transform its left side. Namely, we write
$\bar y -
(\hat\gamma/\gamma)^{-k}
\eta_k = \bar y +
(\hat\gamma/\gamma)^{-k}[\eta_k^0 +
\eta_k^1(\bar x,\mu) + \eta_k^2(\bar y,\mu) +
\eta_k^3(\bar x,\bar y,\mu)]$ where
$\eta_k^1(0,\mu)=0,\eta_k^2(0,\mu)=0$ and
$\eta_k^3 = O(\|\bar x\bar y\|)$. Next, we transfer  constant term
$(\hat\gamma/\gamma)^{-k}\eta_k^0$ into the right side and join it to $M_k$;
we substitute $\bar x$ from the first two equations of (\ref{tk-1+}) into function $\eta_k^1(\bar x,\mu)$ 
and transfer the obtained function to the right side.
After this, all coefficients
(in the third equation) get additions of order $O(\hat\gamma^{-k})$ and
a new linear term $p_ky = O([\hat\gamma/\gamma]^{-k}) y$ appears.
By the shift of coordinates of the form
$(x,y)\mapsto (x,y) + O([\hat\gamma/\gamma]^{-k})$, we vanish both this linear
term and constant terms in the right sides of the first and second equations.
Note also, that now the left side of the third equation can be written as follows
$$
\begin{array}{l}
\bar y + (\hat\gamma/\gamma)^{-k}O(\bar y) +
(\hat\gamma/\gamma)^{-k}O(\|\bar x\bar y\|) =  \\
\qquad = \bar y(1 + q_k) +
(\hat\gamma/\gamma)^{-k}O(|\bar y|^2) +
(\hat\gamma/\gamma)^{-k}O(\|\bar x\bar y\|)
\end{array}
$$
where $q_k = O([\hat\gamma/\gamma]^{-k})$. After this, we can write
system (\ref{tk-1+}) in the form
\begin{equation}
\label{tk-1++}
\begin{array}{l}
\bar x_1 \; = \; a_{11}\lambda_1^{k} x_{1} +
b_1 y + (\hat\lambda^{k}+ \lambda_1^{2k})O(\|x\|) + \lambda_1^{k}O(|y|) + O(y^2) \;, \\
\bar x_2  \; = \;
a_{21}\lambda_1^{k}x_{1} +  b_2 y +
 (\hat\lambda^{k}+ \lambda_1^{2k})O(\|x\|) + \lambda_1^{k}O(|y|) + O(y^2) \;,    \\
\\
\bar y(1 + q_k) +
(\hat\gamma/\gamma)^{-k}O(|\bar y|^2) +
(\hat\gamma/\gamma)^{-k}O(\|\bar x\bar y\|)
\; = \;
M_k  + d \gamma^{k}(1+s_k)y^2 + \\
\qquad + \lambda_1^{k}\gamma^{k}\left[c_{1} x_{1} +
(|\lambda_1|^{k} + |\hat\lambda/\lambda_1|^{k} + \hat\gamma^{-k})O(\|x\|)\right]
+ \lambda_1^k \gamma^{k}O(\|xy\|)
+
\gamma^{k} O(y^3)\;, \\
\end{array}
\end{equation}
where $s_k = O(|\lambda_1|^k + |\hat\gamma/\gamma|^{-k})$ and
new $M_k$ satisfies (\ref{Mk1}).

\subsection{Proof of item 1 of Lemma~\ref{th3}.}\label{prt3-1}

In the Case I we have $b_1(0)=0, c_1\neq 0$ and $b_2\neq 0$.  We choose $\mu_1$ and
$\mu_2\equiv b_1(\mu)$ as the governing parameters. Consider map  (\ref{tk-1++})
 and introduce new coordinates
$$
\begin{array}{l}
\displaystyle
 x_{1new} = x_{1} +
\frac{1}{c_{1}}
(|\lambda_1|^{k} + |\hat\lambda/\lambda_1|^{k} + \hat\gamma^{-k})O(\|x\|)\;, \;
x_{2new} = x_{2}\;,\; y_{new} = y \;,
\end{array}
$$
i.e. we take as $x_{1new}$ the expression from the square brackets
in the third equation of (\ref{tk-1++}). Then (\ref{tk-1++}) is
rewritten in the form
\begin{equation}
\label{tk-1+-}
\begin{array}{l}
\bar x_1 \; = \; (a_{11} +\hat a_k^1)\lambda_1^{k} x_{1} +
\hat a_{k}x_{2} +
(\mu_2 + \rho_k^1) y  + (\lambda_1^{2k}+ \hat\lambda^{k})O(\|x\|^2) +  O(y^2)  \;, \\

\bar x_2  \; = \; a_{21}\lambda_1^{k}x_{1} + (b_2 + \hat r_k)y +
(\lambda_1^{2k}+ \hat\lambda^{k})O(\|x\|) +  O(y^2) \;,    \\
\\
\bar y(1 + q_k) +
\lambda_1^{k}O(|\bar y|^2) +
\lambda_1^{k}O(\|\bar x\bar y\|)
\; = \;
M_k  + \lambda_1^{k}\gamma^{k}c_{1} x_{1} +
d\gamma^{k}(1+ s_k)y^2 + \\
\qquad + \lambda_1^k\gamma^{k}O(\|xy\|) +
\gamma^{k}O(y^3) \;,
\end{array}
\end{equation}
where $\rho_k^1, \hat r_k = O(\lambda_1^{k})$ and
$\hat a_k = O(\lambda_1^{2k} + \hat\lambda^k)$.
Now we rescale the coordinates as follows
$$
\displaystyle y  = -\frac{\gamma^{-k}(1+ q_k)}{d(1+ s_k)}\;Y \;,\; x_1  = \frac{\gamma^{-k}(1+
q_k)}{c_1d(1+ s_k)(\lambda_1\gamma)^k}\;X_1\;,\; x_2  = - \frac{(b_2 + \hat r_k)(1+
q_k)\gamma^{-k}}{d(1+ s_k)}\;X_2.
$$
Then system (\ref{tk-1+-}) is rewritten in the new coordinates as
follows
\begin{equation}
\label{tkres1}
\begin{array}{l}
\displaystyle \bar X_1 \; = \; - J_k X_2 +
M_2 Y  + O(\lambda_1^k)\;,\; \\
\bar X_2  \; = \; Y + O(\gamma^{-k})\;,\; \\
\bar Y +
 = M_1 - X_1 - Y^2 + O(\gamma^{-k})\;,
\end{array}
\end{equation}
where formula (\ref{M123I}) is valid for $M_1$, $M_2$ and $J_k$. Note, that the coefficient $J_k$
in (\ref{tkres1}) is nothing other as the Jacobian of map (\ref{tkres1}) in the point $X=0,Y=0$,
and, hence, $J_k$ coincides, in the main order in $k$, with the Jacobian of map $T_1T_0^k$.

We introduce the new $X_1$-coordinate as follows: $X_{1new}= X_1 - M_2 X_2$. After this the
rescaled map (\ref{tkres1}) takes the form  (\ref{HI}) with $M_3=J_k$.

\subsection{Proof of item 2 of Lemma~\ref{th3}.}\label{prt3-2}

In the Case II we have $c_1(0)=0$ and $b_1\neq 0$.  We choose $\mu_1$ and $\mu_2\equiv
c_1(\mu)$ as the governing parameters. Consider map  (\ref{tk-1++})
 and introduce the new coordinates as $ x_{1new} = x_{1} \;,\; x_{2new} = x_{2} - (b_2/b_1)
x_1 \;,\; y_{new} = y$. Then (\ref{tk-1++}) recasts as
\begin{equation}
\label{tk-1+n}
\begin{array}{l}
\bar x_1 \; = \; a_{11}\lambda_1^{k} x_{1} +
b_1 y + (\hat\lambda^{k} + \lambda_1^{2k})O(\|x\|) + O(y^2) \;,  \\
\bar x_2  \; = \;
a_{21}^\prime\lambda_1^{k}x_{1} +
(\hat\lambda^{k} + \lambda_1^{2k})O(\|x\|) + O(y^2)    \;, \\
%
\bar y(1 + q_k) +
\lambda_1^{k}O(|\bar y|^2) +
\lambda_1^{k}O(\|\bar x\bar y\|)
\; = \;
M_k + d\gamma^{k}(1+ s_k)y^2 + \lambda_1^{k}\gamma^{k}(\mu_2 + \rho_k^2) x_{1} + \\
\qquad + \hat\lambda^k\gamma^k O(|x_2| + \|x\|^2) +
\lambda_1^k\gamma^k O(\|xy\|) + \gamma^{k} O(y^3) \;,
\\
\end{array}
\end{equation}
where $\rho_k^2 = O(|\lambda_1|^{k} +
|\hat\lambda/\lambda_1|^{k})$, $a_{21}^\prime = a_{21} -
(b_2/b_1)a_{11}$ and $M_k = \gamma^{k}[\mu_1 - \gamma^{-k}y^- +
\lambda_1^{k}\mu_2 x_1^+ +  O(\hat\lambda^k + \lambda_1^{2k})]$ in
the case under consideration. Since $c_1(0)=0$, condition $J_1\neq
0$ (see (\ref{t1-3})) implies that $a_{21}^\prime(0) \neq 0$ and therefore $a_{21}^\prime(\mu) \neq 0$ for small $\mu$.
Now we
rescale the coordinates as follows
$$
\displaystyle  y = -\frac{\gamma^{-k}(1+ q_k)}{d(1+ s_k)}\;Y \;,\; x_1  = - \frac{b_1\gamma^{-k}(1+
q_k)}{d(1+ s_k)}\;X_1 \;,\; x_2  = -\lambda_1^k \gamma^{-k}\frac{b_1a_{21}^\prime(1+ q_k)}{d(1+
s_k)}\;X_2 \;.
$$
After this, we can rewrite (\ref{tk-1+n}) in the following form
\begin{equation}
\label{tkres21}
\begin{array}{l}
\displaystyle \bar X_1 \; = \;  Y + O(\lambda_1^k)\;, \\
\bar X_2  \; = \; \displaystyle X_{1} + O\left( (\lambda_1\gamma)^{-k}\right) \;,\\
\displaystyle \bar Y = M_1 + M_2 X_1 + J_k X_2 - Y^2 + O(\lambda_1^k).
\\
\end{array}
\end{equation}
where formula (\ref{M123II}) is valid for $M_1$, $M_2$ and $J_k$. It completes the proof.

\section*{Aknowledgements}
The paper was supported by grant 14-41-00044 of the RSF.
The first and the second authors were partially supported by grants of RFBR No.13-01-00589, 13-01-97028--povolzhje and 14-01-00344,
the second author was also supported by the Leverhulme Trust grant RPG-279 and the EPSRC Mathematics Platform grant EP/I019111/1,
and the third author was supported by the MEC grant MTM2009-09723 (Spain) and the CIRIT grant
2009 SGR 67 (Spain).


\begin{thebibliography}{30}
%
\bibitem{GMO06} S.V. Gonchenko, J.D. Meiss, I.I. Ovsyannikov. Chaotic dynamics of
three-dimensional H\'enon maps that originate from a homoclinic
bifurcation. \textit{Regul. Chaotic Dyn.}, 2006, vol. 11, pp. 191--212.
%
\bibitem{TS08} D.V. Turaev, L.P. Shilnikov. Pseudo-hyperbolisity and the problem
on periodic perturbations of Lorenz-like attractors. \textit{Russian Dokl. Math.}, 2008, vol. 467, pp. 23--27.
English transl.: \textit{Doklady Mathematics}, 2008, vol. 77, no. 1, pp. 17--22.
%
\bibitem{GOST05} Gonchenko S.V., Ovsyannikov I.I., Sim\'o C. and Turaev D.,
Three-dimensional H\'enon-like maps and wild Lorenz-like attractors. \textit{Int. J. of Bifurc. and Chaos}, 2005, 
vol. 15, pp. 3493--3508.
%
\bibitem{SST93}
A.L. Shilnikov, L.P. Shilnikov, D.V. Turaev. Normal forms and
Lorenz attractors. \textit{Int. J. of Bifurc. and Chaos}, 1993, vol. 3, pp. 1123--1139.
%
\bibitem{ASh86}
A.L. Shilnikov. Bifurcation and chaos in the Morioka-Shimizu system. \textit{Methods of
qualitative theory of differential equations}, Gorky, 1986, pp. 180--193. English translation \textit{Selecta
Math. Soviet.}, 1991, vol. 10, pp. 105--117]; \textit{II}. \textit{Methods of Qualitative Theory and Theory of
Bifurcations}, Gorky, 1989, pp. 130--138.
%
\bibitem{ASh93}
A.L. Shilnikov. On bifurcations of the Lorenz attractor in the
Shimuizu-Morioka model. \textit{Physica D}, 1993, vol. 62, pp. 338--346.
%
\bibitem{GGOT13} S.V. Gonchenko, A.S. Gonchenko, I.I. Ovsyannikov, D.V. Turaev
Examples of Lorenz-like attractors in H\'enon-like maps, 
\textit{Mat. Model. of Nat. Phenom}, 2013, vol. 8, no. 5, pp. 48--70.
%
\bibitem{G13} A.S. Gonchenko. On Lorenz-like attractors in model of a Celtic stone.
\textit{Vestnik UdSU, Math., Mech. and Comp. Sci.}, 2013, vol. 2, pp. 3--11.
%
\bibitem{GG13} A.S. Gonchenko, S.V. Gonchenko. On the existence of Lorenz-like attractors in the
nonholonomic model of a ``Celtic stone''. \textit{Rus. J. of Nonlin. Dyn.}, 2013, vol. 9, no. 1, pp. 77--89.
%
\bibitem{GGK13} A.S. Gonchenko, S.V. Gonchenko, A.O. Kazakov. Richness of Chaotic Dynamics in Nonholonomic Models
of a Celtic Stone. \textit{Regular and Chaotic Dynamics}, 2013, vol. 18, no. 5, pp. 521--538.
%
\bibitem{GGS12} A.S. Gonchenko, S.V. Gonchenko, L.P. Shilnikov. 
Towards scenarios of chaos appearance in three-dimensional maps. 
\textit{Rus. J. Nonlinear Dynamics}, 2012, vol. 8, pp. 3--28.
%
\bibitem{GGKT14}
Gonchenko A.S., Gonchenko S.V., Kazakov A.O., Turaev D.V. The simplest scenarios of onset of chaos in three-dimensional maps.
\textit{Int. J. of Bifurc. and Chaos}, 2005, vol. 15, pp. 3493--3508.
%
%
\bibitem{TS98}
D.V.Turaev, L.P.Shilnikov. An example of a wild strange attractor. \textit{Sb. Math}, 1998, vol. 189, no. 2, pp. 137--160.
%
\bibitem{N79} S.E. Newhouse. 
The abundance of wild hyperbolic sets and non-smooth stable sets for diffeomorphisms.
\textit{IHES Publ. Math.}, 1979, vol. 50, pp. 101--151.
%
\bibitem{ASh83a} V.S. Aframovich, L.P. Shilnikov. Strange attractors and quasiattractors. 
\textit{Nonlinear Dynamics and Turbulence}, eds. G.I.Barenblatt, G.Iooss, D.D.Joseph, Boston, Pitmen, 1983.
%
\bibitem{GST93c}
S.V.Gonchenko, L.P.Shilnikov, D.V.Turaev. Dynamical phenomena
in systems with structurally unstable Poincare homoclinic orbits.
\textit{Russian Acad. Sci. Dokl. Math}, 1993, vol. 47, no. 3, pp. 410--415.
%
\bibitem{GST09}
S.V.Gonchenko, L.Shilnikov, D.Turaev. 
On global bifurcations in three-dimensional diffeomorphisms leading to wild Lorenz-like attractors.
\textit{Regul. and Chaotic Dyn.}, 2009, vol. 14, pp. 137--147.
%
\bibitem{GO10} S.V. Gonchenko, I.I. Ovsyannikov. 
On bifurcations of three-dimensional diffeomorphisms with a non-transversal heteroclinic cycle containing saddle-foci.
\textit{Rus. J. Nonlinear Dynamics}, 2010, vol. 6, pp. 61--77.
%
\bibitem{GO13} S.V. Gonchenko, I.I. Ovsyannikov. On Global Bifurcations of Three-dimensional Diffeomorphisms
Leading to Lorenz-like Attractors. \textit{Mat. Model. of Nat. Phenom}, 2013, vol. 8, no. 5, pp. 71--83.
%
\bibitem{T96}
D.V.Turaev. On dimension of nonlocal bifurcational problems.
\textit{Int. J. of Bifurcation and Chaos}, 1996, vol. 6, no. 5, pp. 919--948.
%
\bibitem{Tat01}
J.C.Tatjer. Three-dimensional dissipative diffeomorphisms
with homoclinic tangencies. \textit{Ergod.Th. and Dynam. Sys.}, 2001, vol. 21,
pp. 249--302.
%
\bibitem{NPT}
Newhouse S.E., Palis J., Takens F. Bifurcations and stability of families of diffeomorphisms.
\textit{Publ. Math. Inst. Haute Etudes Scientifiques}, 1983, iss. 57, pp. 5--71.
%
\bibitem{GGT07}
S.V.Gonchenko, V.S.Gonchenko, J.C.Tatjer. Bifurcations of three-dimensional  diffeomorphisms
with non-simple quadratic homoclinic tangencies and generalized Henon maps. \textit{Regular and Chaotic
Dynamics}, 2007, vol. 12, no. 3, pp. 233--266.
%
\bibitem{GST96}
S.V.Gonchenko, L.P.Shilnikov, D.V.Turaev. Dynamical phenomena
in systems with structurally unstable Poincare homoclinic orbits.
\textit{Interdisc. J. CHAOS}, 1996, vol. 6, no. 1, pp. 15--31.
%
\bibitem{GST08}
S.V.Gonchenko, L.P.Shilnikov, D.V.Turaev. On dynamical
properties of multidimensional diffeomorphisms from Newhouse
regions. \textit{Nonlinearity}, 2008, vol. 21, pp. 923--972.
%
\bibitem{Mira96}
C.Mira, L.Gardini, A.Barugola, J.C.Cathala. \textit{Chaotic dynamics
in two-dimensional noninvertible maps}. World Scientific,
Singapore. 1996.
%
\bibitem{PT1}
A.Pumari\~no, J.C.Tatjer. Dynamics near homoclinic
bifurcations of three-dimensional dissipative diffeomorphisms.
\textit{Nonlinearity}, 2006, vol. 19, pp. 2833--2852.
%
\bibitem{PT2}
A.Pumari\~no, J.C.Tatjer. Attractors for return maps near
homoclinic tangencies of three-dimensional dissipative
diffeomorphisms. \textit{Discrete and Continuous Dynamical Systems, series B}, 2007, vol. 8, no. 4, pp. 971--1006.
%
\bibitem{Arneodo79}
A. Arneodo, P. Coullet, C. Tresser.
Possible new strange attractors with a spiral structure.
\textit{Commun. Math. Phys.}, 1981, vol. 79, iss. 4, pp. 573--579.
%
\bibitem{book}
L.P.Shilnikov, A.L.Shilnikov, D.V.Turaev and L.O.Chua.
\textit{Methods of Qualitative Theory in Nonlinear Dynamics, Part~I}. World
Scientific. 1998.
%
\bibitem{GS90}
S.V.Gonchenko, L.P.Shilnikov. Invariants of
$\Omega$-conjugacy of diffeomorphisms with a nontransversal
homoclinic orbit. \textit{Ukr. Math. J}, 1990, vol. 42, no. 2, pp. 134--140.
%
\bibitem{GS92}
S.V.Gonchenko, L.P.Shilnikov. On moduli of systems with a
nontransversal Poincare homoclinic orbit. \textit{Russian
Acad. Sci. Izv. Math}, 1993, vol. 41, no. 3, pp. 417--445.
%
\bibitem{HPS}
M.W.Hirsch, C.C.Pugh, M.Shub. \textit{Invariant manifolds. Lecture
Notes in Math}, 1977, vol. 583, Springer-Verlag, Berlin.
%

\end{thebibliography}
\end{document}